\documentclass{article}

\title{\LARGE \textbf{Simple Proofs of two Dirac-type Theorems Involving Connectivity}}
\author{Carlen Mosesyan, Mher Nikoghosyan, Zhora Nikoghosyan}

\begin{document}

\maketitle

\begin{abstract}
Simple and shorter proofs of two Dirac-type theorems involving connectivity are presented.   \\

\end{abstract}

We consider only finite undirected graphs without loops or multiple edges. Let $G$ be a graph of order $n$ with minimum degree $\delta$, connectivity $\kappa$ and circumference $c$ - the length of a longest cycle $C$ in $G$. Then $C$ is a Hamilton cycle if $|C|=n$ and is a dominating cycle if $G\backslash C$ is edgeless.   

In 1981, two Dirac-type theorems appeared involving connectivity $\kappa$.\\ 

\noindent\textbf{Theorem 1} [3]. In every 3-connected graph, $c\geq\min \{n,3\delta-\kappa\}$.\\

\noindent\textbf{Theorem 2} [3]. Every 2-connected graph with $\delta\geq (n+\kappa)/3$ is hamiltonian.\\

A short proof of Theorem 2 was given in [1] based on maximal independent sets. Some other proofs of Theorems 1 and 2 are given in terms of degree sums. 

In this note we present much more simple and shorter proofs of theorems 1 and 2 based mainly on standard arguments and the following two theorems.\\

\noindent\textbf{Theorem 3} [4]. Let $G$ be a 3-connected graph. Then either $c\geq 3\delta-3$ or every longest cycle in $G$ is a dominating cycle.\\

\noindent\textbf{Theorem 4} [2]. Let $G$ be a 2-connected graph with $\delta\geq (n+2)/3$. Then every longest cycle in $G$ is a dominating cycle.\\

The set of vertices of a graph $G$ is denoted by $V(G)$ and the set of edges by $E(G)$. For $S$ a subset of $V(G)$, we denote by $G\backslash S$ the maximum subgraph of $G$ with vertex set $V(G)\backslash S$. For a subgraph $H$ of $G$ we use $G\backslash H$ short for $G\backslash V(H)$. We denote by $N(x)$ the neighborhood of a vertex x in a graph $G$. We write a cycle $C$ of $G$ with a given orientation by $\overrightarrow{C}$. For $x,y\in V(C)$, we denote by $x\overrightarrow{C}y$ the subpath of $C$ in the chosen direction from $x$ to $y$. For $x\in V(C)$, we denote the successor of $x$ on $\overrightarrow{C}$ by $x^+$. For $X\subset V(C)$, we define $X^+=\{x^+|x\in X\}$.

\noindent\textbf{Lemma 1}. Let $G$ be a graph and $S$ a minimum cut-set in $G$. If every longest cycle in $G$ is a dominating cycle, then either $c\geq 3\delta-\kappa+1$ or there exist a longest cycle $C$ with $S\subseteq V(C)$.\\

\noindent\textbf{Proof}. Choose a longest cycle $C$ in $G$ such that $|V(C)\cap S|$ is as great as possible and assume that $S\not\subseteq V(C)$ with $x\in S\backslash V(C)$. Since $C$ is dominating, $N(x)\subseteq V(C)$.       Let $\xi_1,...,\xi_t$ be the elements of $N(x)$, occuring on $\overrightarrow{C}$ in a consecutive order and dividing $C$ into segments $I_i=\xi_i\overrightarrow{C}\xi_{i+1}$ $(i=1,...,t)$, where $\xi_{t+1}=\xi_1$. Put $M_1=\{\xi_i|V(I_i)\cap S\subseteq\{\xi_i,\xi_{i+1}\}\}$ and $M_2=N(x)\backslash M_1$.      Clearly $|M_2|\leq \kappa-1$. Since $C$ is extreme, $N(x)$, $N^+(x)$ and $M_1^{++}$ are pairwise disjoint in $V(C)$. Then the result follows from $|M_1^{++}|=|N(x)|-|M_2|\geq\delta-\kappa+1$.
  \qquad $\Delta$\\

\noindent\textbf{Proof of Theorem 1}. Let $G$ be a 3-connected graph, $S$ be a minimum cut-set in $G$ and let $H_1,...,H_h$ be the components of $G\backslash S$. The result holds immediately if $c\geq3\delta-3$, since $3\delta-3\geq3\delta-\kappa$. Otherwise, by Theorem 3, every longest cycle in $G$ is a dominating cycle. Let $C$ be any one with $x\in V(G\backslash C)$. By Lemma 1, $S\subseteq V(C)$ and we can assume w.l.o.g. that $x\in V(H_1)$. Put $W_1=N(x)\cup N^+(x)$. Clearly $|W_1|\geq2\delta$ and it remains to find $W_2\subseteq V(C)$ such that $W_1\cap W_2=\emptyset$ and $|W_2|\geq\delta-\kappa$. Suppose first that $W_1\subseteq V(H_1)\cup S$. If $V(H_2)\subseteq V(C)$, then take $W_2=V(H_2)$. Otherwise, there exist $y\in V(H_2\backslash C)$ with $N(y)\subseteq V(C)$ and we can take $W_2=N(y)\backslash S$. Now let $W_1\not\subseteq V(H_1)\cup S$ and choose $z\in N^+(x)\cap V(H_2)$. If $N(z)\subseteq V(C)$, then take $W_2=N(z)\backslash S$, since $N(z)\cap N^+(x)=\emptyset$    (by standard arguments). Otherwise, choose $w\in N(z)\backslash V(C)$. Clearly $N(w)\subseteq V(C)$, $w\in V(H_2)$ and $N(z)\cap N^+(x)\subseteq \{z\}$. Then by taking $W_2=(N(w)\backslash \{z\})\backslash (S\backslash \{z^-\})$ we complete the proof.   \qquad $\Delta$\\

\noindent\textbf{Proof of Theorem 2}. Assume the converse. Let $G$ be a non hamiltonian 2-connected graph with $\delta\geq(n+\kappa)/3$ and $S$ be a minimum cut-set in $G$. Since $\delta\geq(n+\kappa)/3\geq(n+2)/3$, by Theorem 4, every longest cycle in $G$ is a dominating cycle. By Lemma 1, $G$ contains a longest cycle $C$ with $S\subseteq V(C)$. As in proof of Theorem 1, $c\geq3\delta-\kappa$, contradicting the fact that $\delta\geq(n+\kappa$)/3.     \qquad $\Delta$

\end{document}